# On a new scale of regularity spaces with applications to Euler's equations

Eitan Tadmor*

**Abstract**

We introduce a new ladder of function spaces which is shown to fill in the gap between the weak $L^{p\infty}$ spaces and the larger Morrey spaces, $M^p$.

Our motivation for introducing these new spaces, denoted $\vee^{pq}$, is to gain a more accurate information on (compact) embeddings of Morrey spaces in appropriate Sobolev spaces. It is here that the secondary parameter $q$ (— and a further logarithmic refinement parameter $\alpha$, denoted $\vee^{pq}(\log\vee)^{\alpha}$) gives a finer scaling, which allows us to make the subtle distinctions necessary for embedding in spaces with a fixed order of smoothness.

We utilize an $H^{-1}$-stability criterion which we have recently introduced in [21], in order to study the strong convergence of approximate Euler solutions. We show how the new refined scale of spaces, $\vee^{pq}(\log\vee)^{\alpha}$, enables us approach the borderline cases which separate between $H^{-1}$-compactness and the phenomena of concentration-cancelation. Expressed in terms of their $\vee^{pq}(\log\vee)^{\alpha}$ bounds, these borderline cases are shown to be intimately related to uniform bounds of the total (Coulomb) energy and the related vorticity configuration.

# Contents

*Department of Mathematics, UCLA, Los-Angeles CA 90095, USA. *Email:* tadmor@math.ucla.edu

# 1 Introduction

We introduce a new ladder of function spaces which is shown to fill in the gap between Marcenkwetiz weak-$L^p$ spaces, and the larger Morrey spaces, $M^p$. The former measure the total mass of measurable $f$'s over *arbitrary sets*; the latter measure the total mass of $f$'s over *arbitrary balls.* The newly introduced scale of spaces, denoted $\vee^{pq}$, measures a $\ell_q$-weighted distribution of the total mass of a measurable $f$ over an arbitrary *collection of disjoint balls*

Our motivation for introducing these new spaces, denoted $\vee^{pq}$, is to gain a more accurate information on (compact) embeddings of Morrey spaces in appropriate Sobolev spaces. It is here that the secondary parameter $q$ (and a further logarithmic refinement parameter $\alpha$, denoted $\vee^{pq}(\log \vee)^\alpha$) give a finer scaling, which allows us to make the subtle distinctions necessary for embedding in spaces with a fixed order of smoothness.

In §2 we prove that the new scale of spaces, $\vee^{pq}$, $q \geq p \geq 1$, interpolates the gap between $\vee^{pp} \subset L^{p\infty}$ and $\vee^{p\infty} = M^p$. The further logarithmic refinement, $\vee^{pq}(\log \vee)^\alpha$ is introduced in order to address the case $p = 1$. In §3 we study the compact embeddings of $\vee^{pq}(\log \vee)^\alpha(I\!R^N)$. The compactness results are accomplished by precise characterization of the decay of the wavelet coefficients for $\vee^{pq}(\log \vee)^\alpha$-functions. We are particularly interested in $H^{-1}$-compactness. For $\Omega \subset I\!R^N$, it is shown that $\vee^{p2}(\log \vee)^\alpha(\Omega)$ is $H^{-1}$-compact if $p > 2N/(N+2)$, or, if $p = 2N/(N+2)$ and $\alpha > 1/2$. This should be compared with the compact imbedding of Morrey spaces, consult [21, Theorem 4.2], which states that $M^p(\Omega) \overset{\text{comp}}{\hookrightarrow} H^{-1}$ for the restricted range of $p > N/2$. Equipped with the new scale of spaces $\vee^{p2}(\log \vee)^\alpha(I\!R^N)$ we are now able to resolve the question of compactness for the gap of $p$'s, $p \leq N/2 < 2$. Specifically, we show that the question of $H^{-1}_{\text{loc}}(I\!R^N)$-compactness is characterized by the borderline cases of $X_2 := \vee^{12}(\log \vee)^{1/2}_c(I\!R^2)$ in the two-dimensional case, and $X_3 := \vee^{\frac{6}{5}2}_c(I\!R^3)$ in the three-dimensional case.

The new scale of spaces, through the precise characterization of their $H^{-1}$-compactness, is put into use in §4, where we discuss approximate solution of the incompressible Euler equations. Recently, we introduced in [21] a sharp local condition for the lack of concentrations in (— and hence the $L^2$ convergence of) sequences of such approximate solutions. Simply stated, the sequence of associated vorticities is required to be $H^{-1}_{\text{loc}}$-compact, and it is in this context that the $\vee^{pq}$-bounds are shown to play a fundamental role. Indeed, in both the $N=2$ and the $N=3$-dimensional cases, we show that the corresponding $X_N$-bounds on the vorticities are intimately related to the uniform bound on the Coulomb energy of the solutions and their vorticity configurations.

In the two-dimensional case we end up with a rather complete classification which is summerize in the following statement (consult Corollary 4.1 and Theorem 4.1 below).

**Theorem**. *Let $\{u^\varepsilon(\cdot,t)\}$ be a family of approximate solutions of the 2D Euler equations, and assume that the corresponding sequence of vorticities $\{\omega^\varepsilon(\cdot,t)\}$ is uniformly bounded in $\widetilde{\vee}^{12}(\log \widetilde{\vee})^\alpha_c(I\!R^2)$, $\alpha > 0$.*

{i} (No concentration). *If $\alpha > 1/2$, then $\{u^\varepsilon(\cdot,t)\}$ is strongly compact in $L^\infty([0,T]; L^2_{loc}(I\!R^2))$, with a strong $L^2$-limit, $u(\cdot,t)$, which is a weak solution of the 2D Euler equations.*

{ii} (Concentration-Cancellation). *If* $\alpha \in (0, \frac{1}{2}]$*, then* $\{u^\varepsilon(\cdot,t)\}$ *has a* $L^2$*-weak limit,* $u(\cdot,t)$ *which is a finite-energy solution of the 2D Euler equations.*

One signed measures, say $\omega^\varepsilon(\cdot,0) \geq 0$ are shown to be bounded in $\widetilde{\vee}^{12}(\log \widetilde{\vee})_c^{1/2}(I\!R^2)$ (consult Lemma 4.1 below), and thus part (ii) of the above theorem recasts an extended version of Delort result [14] in the language of $\vee$-spaces.

The new ladder of spaces establishes a direct linkage between questions related to the configuration of the $N$-dimensional pseudo-energy and regularity in the borderline case $X_N$. The configuration of three-dimensional vorticity, $\boldsymbol{\omega}^\varepsilon(\cdot,t)$ involves local stretching effects and nonlinear energy saturation associated with small sets of increasingly intense vorticity. The following result relates these issues to the uniform bound in the borderline case $X_3 = \vee^{\frac{6}{5}2}$.

**Theorem**. *Let* $\{u^\varepsilon(\cdot,t)\}$ *be a family of approximate solutions of the 3D Euler equations and assume that the corresponding sequence of compactly supported vorticities,* $\{\boldsymbol{\omega}^\varepsilon(\cdot,t)\}$ *satisfies the local alignment condition (4.31). Then the following holds,*

$$\|\boldsymbol{\omega}^\varepsilon(\cdot,t)\|_{\vee^{\frac{6}{5}2}(\Omega)} \leq Const., \qquad \Omega \subset I\!R^3.$$

We close by noting that the new scale of spaces, $\widetilde{\vee}^{pq}(\Omega)$, is not necessarily comparable with the scale of Morrey spaces, $\widetilde{M}^r(\Omega)$, unless additional information, e.g., the packing measure of $\Omega$ is provided. Thus for example, $\widetilde{M}^{3/2}(\Omega)$ is borderline Morrey space for $H^{-1}_{\rm loc}(I\!R^3)$ compactness, which was shown by Giga & Miyakawa [15] to guarantee the existence of the related 3D Navier-Stokes solutions. Compared with the corresponding borderline case $X_3(\Omega) = \widetilde{\vee}^{\frac{6}{5}2}(\Omega)$, we find in Corollary 3.1 below that the latter is larger, $\widetilde{M}^{3/2}(\Omega) \subset X_3(\Omega)$, for $\Omega$'s with finite packing measure so that $\pi^{h_1}(\Omega) < \infty$.

**Acknowledgment** I am grateful to Ron DeVore for many fruitful discussions that led to this work and for his critical reading of the manuscript.
Research was supported in part by ONR grant N00014-91-J-1076 and NSF grant #DMS97-06827.

# 2 The spaces $\vee^{pq}(\log \vee)^\alpha(\Omega)$

Given a domain $\Omega \subset I\!R^N$, we consider the set $\mathcal{B}(\Omega)$ of all collections of mutually disjoint balls contained in $\Omega$, $\mathcal{B}(\Omega) = \{B_j \mid \cup B_j \subseteq \Omega\}$, balls with sufficiently small radius $B_j = B_{R_j}(x_j),\ R_j \leq R_0 < 1/2$.

**Definition 2.1** *The space* $\vee^{pq}(\Omega),\ 1 \leq p \leq q \leq \infty$*, consists of all* $f$*'s in* $L^1_{loc}(\Omega)$ *such that for all collections,* $\{B_j\} \subset \mathcal{B}(\Omega)$*, the following estimate holds*

$$\sup_{\{B_j\}\subset\mathcal{B}(\Omega)} \left( \sum_j \left( R_j^{-N/p'} \int_{B_j} |f(x)|dx \right)^q \right)^{1/q} \leq Const, \quad 1 \leq p \leq q \leq \infty. \qquad (2.1)$$

The smallest of such *Const*'s in (2.1) is the $\vee^{pq}$-norm of $f$. Thus, if let $⨍_{B_R(x_0)}|f(x)|dx$ denote the average mass of $|f|$ over the ball $B_R(x_0)$ centered at $x_0$, and $\bar{f}=(\bar{f}_1,\bar{f}_2,\ldots,)$ denote the vector of averages, $\bar{f}_j := ⨍_{x\in B_j}|f(x)|dx$, then

$$\|f\|_{\vee^{pq}(\Omega)} := \sup_{R_j<R_0}\|\{R_j^{N/p}\bar{f}_j\}\|_{\ell^q}, \quad q\geq p. \tag{2.2}$$

Occasionally, we shall need a further logarithmic refinement

$$\vee^{pq}(\log\vee)^{\alpha}(\Omega) := \{f\in L^1(\Omega)\ |\ \sup_{R_j<R_0}\|\{R_j^{N/p}|\log R_j|^{\alpha}\bar{f}_j\}\|_{\ell^q}<\infty\}, \quad q>p. \tag{2.3}$$

We abbreviate $\vee^{pq,\alpha}=\vee^{pq}(\log\vee)^{\alpha}$. We shall also need the corresponding extension dealing with bounded measures, $\mu\in\mathcal{BM}(\Omega)$. With $\bar{\mu}_j := |\mu|(B_j)/|B_j|$ we set

$$\widetilde{\vee}^{pq,\alpha} = \{\mu\in\mathcal{BM}(\Omega)\ |\ \sup_{R_j<R_0}\|\{R_j^{N/p}|\log R_j|^{\alpha}\bar{\mu}_j\}\|_{\ell^q}, \quad q\geq p.$$

For $\Omega=I\!R^N$, the space $\vee^{pq,\alpha}_{\rm loc}$ is defined as the Fréchet space determined by the family $\{\|f\|_{\vee^{pq,\alpha}(B_k(0))}\}_{k\in I\!N}$.

The norm $\|f\|_{\vee^{pq,\alpha}(\Omega)}$ quantifies the (ir-)regularity of $f$ by measuring a weighted distribution of its singularities, distributed over a 'packing' of $\Omega$ by a covering of balls. Clearly, the use of balls in these definitions is not essential, and they can be replaced, for example by sequences of non-thin cubes, $\{\mathcal{C}_j\}$, for which $(diam\ \mathcal{C}_j)^N\leq Const|\mathcal{C}_j|$. Thus, if a bounded $\Omega\subset I\!R^N$ is covered by a lattice of disjoint cubes, $\Omega\subset\cup\mathcal{C}_j,\ \mathcal{C}_j(\cdot)=\mathcal{C}(\cdot+j),\ j\in Z^N$, each of which of size $|\mathcal{C}_j|=R^N$, then $f\in\vee^{pq,\alpha}(\Omega)$ implies

$$\left(\sum_j\Big(\int_{\mathcal{C}_j}|f(x)|dx\Big)^q\right)^{1/q}\leq R^{N/p'}|\log R|^{-\alpha}, \quad \mathcal{C}_j(\cdot):=\mathcal{C}(\cdot+j),\ j\in Z^N. \tag{2.4}$$

We note in passing that as we refine the covering, say by a dyadic refinment of the lattice in (2.4), the corresponding $\vee$-sum, $\|\{R_j^{N/p}\bar{f}_j\}\|_{\ell^q}$ need not increase for $q\geq p$.

We want to place the scale of new spaces, $\vee^{pq,\alpha}$, in the context other known spaces, and this is carried out in terms of the known Lorentz-Zygmund and Morrey spaces. Here is a brief readers' digest which will enable us to introduce the necessary notations, and we refer the redear to [2] for a detailed description.

Let $f^*$ denote the usual decreasing rearrangement of $f$. For a bounded $\Omega\subset I\!R^n$, the space $L^{pq,\alpha}(\Omega)=L^{pq}(\log L)^{\alpha}(\Omega)$ consists of all measurable functions $f$'s such that $\Big(\int_{s=0}^{|\Omega|}[s^{1/p}|\log s|^{\alpha}f^*(s)]^q ds/s\Big)^{1/q}<\infty$; we shall be exclusively concerned with the weak spaces corresponding to $q=\infty$, where

$$\|f\|_{L^{p\infty,\alpha}(\Omega)} = \sup_{s\leq|\Omega|} s^{1/p}|\log s|^{\alpha}f^*(s).$$

For consistency of notations, however, we retain here the secondary index $q = \infty$, and we refer the interested reader to [1],[2], for a detailed study of the logarithmic refinement indexed with $\alpha > 0$.

If we replace $f^*$ with its maximal function, $f^{**} := \frac{1}{s}\int_0^s f^*(r)dr$, we obtain the closely related Lorentz-Zygmund spaces $L^{(pq,\alpha)}(\Omega) = \{f \mid \|s^{1/p}|\log s|^\alpha f^{**}\|_{L^q(ds/s)} < \infty\}$. The $L^{(pq,\alpha)}$'s are rearrangement-invariant spaces which include as special cases both the Lorentz spaces, $L^{(pq)} = L^{(pq,0)}$, and the logarithmic Orlicz spaces $L^{(11,\alpha-1)} = L(\log L)^\alpha(\Omega)$ [1, Theorem 11.1], Again, we are exclusively interested here in the case of secondary index $q = \infty$: using the maximality of $f^{**}(s) = \sup_{|E|=s} \fint_E |f|$, we find that $L^{(p\infty,\alpha)}$ consists of all $f$'s such that

$$L^{(p\infty,\alpha)}(\Omega) = \{f \mid \int_E |f(y)|dy \leq Const \cdot |E|^{1/p'} |\log |E|\,|^{-\alpha}, \quad \forall E \subset \Omega, |E| < E_0 < 1\}. \tag{2.5}$$

We note in passing that $L^{(pq,\alpha)}$ coincide with $L^{pq,\alpha}$ for $p > 1$, [1, Corollary 8.2]. For $p = 1$, however, the spaces $L^{(1q,\alpha)}$ ( – denoted $\mathcal{L}^{1q}(\log \mathcal{L})^\alpha$ in [1, §11]) are strictly smaller than the corresponding $L^{1q,\alpha}$. Thus, with $\alpha = 0$ for example, the $L^{(1q)}$'s are varying between $L^{(11)} = L(\log L)$ and $L^{(1\infty)} = L^1$.

Finally, if we replace in (2.5) the arbitrary sets $E$'s by balls, we enlarge the Lorentz-Zygmund spaces, arriving the scale of Morrey spaces,

$$M^{p,\alpha}(\Omega) := \{f \in L^1(\Omega) \mid \int_{B_R(x_0)\subset\Omega} |f(x)|dx \leq CR^{N/p'}|\log R|^{-\alpha}, \quad \forall R \leq R_0 < 1\}. \tag{2.6}$$

The case $\alpha = 0$ yields the classical Morrey space $M^p$, e.g., [16], [15]; the logarithmic refinement of $M^{p,\alpha}$ was put into a recent use in [21], motivated by the corresponding logarithmic refinement in Lorentz-Zygmund spaces. Following [15] and [21], we let $\widetilde{M}^{p,\alpha}(\Omega)$ denote the corresponding Morrey scale for bounded measures, $\mu \in \mathcal{BM}(\Omega)$

$$\|\mu\|_{\widetilde{M}^{p,\alpha}} := \sup_{R<R_0<1} \left[R^{-N/p'}|\log R|^\alpha |\mu|(B_R(x))\right].$$

We now turn to discuss the scale of spaces $\vee^{pq}$ for $p \leq q \leq \infty$. Their definition in (2.1) makes apparent the role of the parameter $q$ as the usual secondary index, so that $\vee^{pq}$ form a 'scale' of intermediate spaces between $\vee^{p\infty}$ and $\vee^{pp}$. Indeed one can interpolate an-$\vee^{pq}$ bound

$$\|f\|_{\vee^{pq}} \leq \|f\|^\theta_{\vee^{pp}} \cdot \|f\|^{1-\theta}_{\vee^{p\infty}}, \quad \theta = p/q \leq 1. \tag{2.7}$$

More precisely, using real interpolation arguments along the lines of e.g., [9, theorem 7.5], one finds $\vee^{pq}$ as an interpolation space of $\vee^{p\infty}$ and $\vee^{pp}$,

$$\vee^{pq} = (\vee^{p\infty}, \vee^{pp})_{\theta,q}, \quad \theta = p/q \leq 1.$$

It is therefore enough to consider the two end cases — $q = p$ and $q = \infty$. We start with the latter.

Clearly, $\vee^{p\infty,\alpha}$ consists of all $L^1_{\rm loc}(\Omega)$ $f$'s whose behavior is determined by their average mass on just *one* ball, i.e., (2.6) holds.

**Lemma 2.1** *For $p \geq 1$ we have*

$$\vee^{p\infty,\alpha}(\Omega) = M^{p,\alpha}(\Omega), \qquad p \geq 1. \tag{2.8}$$

Next we turn to discuss the spaces $\vee^{pp,\alpha}$, which are shown to be in between the Lorentz-Zygmund spaces $L^{pp,\alpha} \equiv L^p(\log L)^\alpha$ (– consult [1, Corllary 10.2]) and $L^{p\infty,\alpha} \equiv L^{p\infty}(\log L)^\alpha$. The following lemma is in the heart of this matter.

**Lemma 2.2** *For $p \geq 1$ we have*

$$L^p(\log L)^\alpha(\Omega) \subset \vee^{pp,\alpha}(\Omega) \subset L^{p\infty}(\log L)^\alpha(\Omega), \qquad p \geq 1, \ \alpha \geq 0. \tag{2.9}$$

**Proof**. Consider an arbitrary open measurable set $E \subseteq \Omega$ of size $|E| = t < 1$. To verify the right half of (2.9), we need to estimate the decay rate of $\int_E |f(x)|dx$ as $t \downarrow 0$. To this end, we cover $E$ by the family of interior balls $\{B_{R_x}(x) \subset E\}$. By Vitali's covering lemma, e.g., [28, §1.6], we can select a subfamily of countably many disjoint balls, $\{B_j = B_{R_j}(x_j) \mid \cup B_j \subset E\}$, which cover at least a fixed fraction of $E$, namely, the complement of $E_1 := \cup_j B_{R_j}(x_j)$ does not exceed $|E - E_1| \leq \theta t$ with $\theta = (4/5)^N$.

We write

$$\int_E |f(x)|dx = \int_{E-E_1} |f(x)|dx + \sum_j \int_{B_j} |f(x)|dx. \tag{2.10}$$

Assuming that $f \in \vee^{pp,\alpha}(\Omega)$, then the last summation on the right does not exceed

$$\begin{aligned}\sum_j \int_{B_j} |f(x)|dx &\leq \left(\sum_j \left(R_j^{-N/p'}|\log R_j|^\alpha \int_{B_j} |f(x)|dx\right)^p\right)^{1/p} \cdot \left(\sum_j |\log R_j|^{-\alpha p'} R_j^{Np'/p'}\right)^{1/p'} \\ &\leq Const \cdot t^{1/p'}|\log t|^{-\alpha}, \qquad Const = N^\alpha \|f\|_{\vee^{pp,\alpha}(\Omega)}.\end{aligned}$$

Next, consider the maximal function $F(t) := \sup_{|E|=t} \int_E |f(x)|dx$. Using the fact that $|E - E_1| \leq \theta t$ together with (2.11), then (2.10) yields

$$F(t) \leq F(\theta t) + Const \cdot t^{1/p'}|\log t|^{-\alpha}. \tag{2.11}$$

Recalling that $F(t)$ is in fact the primitive of the decreasing rearrangement $f^*$, $F(t) = \int_0^t f^*(s)ds$, the desired $\|f\|_{L^{p\infty}}$-bound follows from (2.11)

$$f^*(t) \leq \frac{F(t) - F(\theta t)}{(1-\theta)t} \leq Const \cdot t^{-1/p}|\log t|^{-\alpha}.$$

For the reversed implication on the left of (2.9), we use Hölder inequality which yields the following straightforward $\vee^{pp}$-bound for $L^p$ functions,

$$\begin{aligned}\sum_j \left(R_j^{-N/p'}\int_{B_j}|f(x)|dx\right)^p &\leq \sum_j R_j^{-Np/p'}\int_{B_j}|f(x)|^p dx \times R_j^{Np/p'} = \\ &= \int_{\cup B_j}|f(x)|^p dx \leq \|f\|^p_{L^p}, \quad p\geq 1, \end{aligned} \tag{2.12}$$

and thus, the LHS of (2.9) with $\alpha=0$ follows. For general $\alpha>0$ we need a logarithmic refinment based on the duality between $L^p(\log L)^\alpha$ and $L^{p'}(\log L)^{-\alpha}$, consult e.g., [2, Corollary 8.15], [1, Theorem 8.4], yielding

$$\begin{aligned}&\sum_j \left(R_j^{-N/p'}|\log R_j|^\alpha\int_{B_j}|f(x)|dx\right)^p \leq \\ &\leq \sum_j R_j^{-Np/p'}|\log R_j|^{\alpha p}\times\|f(x)\|^p_{L^p(\log L)^\alpha(B_j)}\cdot\left(\int_{t=0}\left[(1+|\log t|)^{-\alpha}1_{0\leq t\leq R_j^N}\right]^{p'}dt\right)^{p/p'} \leq \\ &\leq Const\sum_j R_j^{-Np/p'}|\log R_j|^{\alpha p}\times\|f(x)\|^p_{L^p(\log L)^\alpha(B_j)}\cdot R_j^{Np/p'}|\log R_j|^{-\alpha p}\leq \\ &\leq Const\|f(x)\|^p_{L^p(\log L)^\alpha(\cup B_j)}, \quad p\geq 1,\ \alpha\geq 0. \end{aligned} \tag{2.13}$$

Thus, the $\vee^{pp,\alpha}$ bound of $f$ implies that the LHS of (2.9) holds. ∎

*Remarks.*

1. We note in passing an alternative derivation of (2.9). Setting $F^{(p,\alpha)}(t) := t^{-1/p'}|\log t|^\alpha F(t)$, then (2.11) yields

$$F^{(p,\alpha)}(t)\leq \theta^{1/p'}\left|\frac{\log t}{\log(\theta t)}\right|^\alpha F^{(p,\alpha)}(\theta t)+Const.$$

Successive application of this recursion relation yields

$$\begin{aligned}F^{(p,\alpha)}(t) &\leq \sum_k^\infty \theta^{k/p'}\left|\frac{\log t}{\log(\theta^k t)}\right|^\alpha \leq \\ &= |\log t|^\alpha\sum_k\frac{\theta^{k/p'}}{(k|\log\theta|+|\log t|)^\alpha}\leq\begin{cases} Const., & p>1 \\ Const\cdot|\log t|, & p=1,\alpha>1\end{cases}.\end{aligned}$$

For $p>1$, we conclude, as before, $\vee^{pp,\alpha}\subset L^{(p\infty,\alpha)}=L^{p\infty,\alpha}$ — the logarithmic refinement corresponding to (2.9). For $p=1,\alpha>1$, however, this approach only yields $\vee^{11,\alpha}\subset L^{(1\infty,\alpha-1)}$, whereas the derivation of Lemma 2.2 led to a tighter bound in terms of $L^{1\infty,\alpha}$. We note that the space $L^{1\infty,\alpha}$ is indeed smaller (at least for $\alpha > 1$) than the space $L^{(1\infty,\alpha-1)}$ [1, Theorem 12.1].

2. The following example, due to R. DeVore, [7], shows that $\vee^{pp}(I\!R_+)$ lies *strictly* inside $L^{p\infty}(I\!R_+)$. To this end observe that the averages of the $L^{p\infty}$ function $f(x) = x^{-1/p}$, averaged over the dyadic intervals $I_j = [2^{-j}, 2^{-j+1}]$, are given by $\bar{f}_j = ⨍_{I_j} |f(x)|dx = c_p 2^{j/p}$, and hence $\{2^{-j/p}\bar{f}_j \equiv Const\} \notin \ell_p$. In fact, this shows that $L^{p\infty} \not\subset \vee^{pq},\ q < \infty$.

3. For a different kind of inclusion relations in terms of Besov spaces we refer to (3.18) below, asserting that $\vee^{pp'}(\Omega) \subset B_\infty(L^{p'}(\Omega))$ for $p \leq 2$.

In summary, we see that the new spaces $\vee^{pq,\alpha}$ offer a new ladder which covers the gap between the weak Lorentz-Zygmund spaces corresponding to $q = p$, and the larger Morrey spaces corresponding to $q = \infty$, namely

$$L^{pp,\alpha} \subset \vee^{pp,\alpha}_{\text{loc}} \subset L^{p\infty,\alpha} \ldots \vee^{pq,\alpha}_{\text{loc}} \ldots \subset \vee^{p\infty,\alpha}_{\text{loc}} = M^{p,\alpha}_{\text{loc}}, \quad p \geq 1. \tag{2.14}$$

# 3 Compact Imbeddings

Our motivation for introducing the new spaces $\vee^{pq,\alpha}$ was to gain a more accurate information on (compact) embeddings of Morrey spaces in appropriate Sobolev spaces. It is here that the secondary parameters $q$ and $\alpha$ give a finer scaling, which allows us to make the subtle distinctions necessary in embedding in spaces with a fixed order of smoothness. To avoid an excessive amount of indices, we begin with a prototype configuration, referring to the specific situation encountered in [21]. The general case will be stated later (in Theorem 3.2 below).

According to [21, Theorems 4.2 & 4.3], the Morrey spaces $\widetilde{M}^{p,\alpha}(\Omega)$ are precompact in $H^{-1}(\Omega)$ as long as

$$\widetilde{M}^{p,\alpha}(\Omega) \stackrel{\text{comp}}{\hookrightarrow} H^{-1}(\Omega), \quad \left(p - \frac{N}{2}\right)_+ + (\alpha - 1)_+ > 0. \tag{3.15}$$

We distinguish between two borderline cases.

• In the two-dimensional case, we find that $\widetilde{M}^{1,\alpha}(I\!R^2)$ is precompact in $H^{-1}(I\!R^2)$ for $\alpha > 1$. On the other hand, counterexamples constructed in [11],[22] show that $\widetilde{M}^{1,1/2}(I\!R^2) \cap \mathcal{BM}^+_c(I\!R^2)$ is not compactly imbedded in $H^{-1}(I\!R^2)$. Thus, the gap $1/2 < \alpha < 1$ remains open with regard to the question of compact embedding of $\widetilde{M}^{1,\alpha}(I\!R^2)$ in $H^{-1}_{\text{loc}}(I\!R^2)$.

• The gap is even wider for $p > 1$. Considering the Lebesgue/Lorentz hierarchy (here we ignore the logarithmic subscaling, taking $\alpha = 0$), one finds the critical Lebesgue exponent $(p^*)' = \frac{2N}{N+2}$, so that all $L^{p,\infty}_c(I\!R^N)$ with $p > \frac{2N}{N+2}$ are compactly imbedded in $H^{-1}_{\text{loc}}(I\!R^N)$. The Morrey hierarchy is different: according to (3.15), Morrey spaces $M^p_c(I\!R^N)$ are $H^{-1}_{\text{loc}}(I\!R^N)$-compact for a smaller range of exponents with $p > N/2$. Though the Morrey spaces are bigger than the corresponding weak-$L^p$, $L^{p,\infty} \subset M^p$, they both admit the same scaling. Thus, for example, with $N = 3$ we are left with the open question with regard to the 'correct' scaling exponent within the intermediate gap $\frac{6}{5} < p < \frac{3}{2}$, which will suffice for compact imbedding in $H^{-1}_{\rm loc}(I\!R^3)$.

Equipped with the new scale of *intermediate* spaces $\vee^{pq,\alpha}$, we are able to address the question of compactness for the above gaps, by sharpening (3.15) as follows.

**Theorem 3.1** *Let $\Omega \subset I\!R^N$ be a bounded domain and let $\{f^\varepsilon\} \subset C_c^\infty(\Omega)$ be a bounded sequence in $\vee^{p2,\alpha}(\Omega)$. If either:*

*(a) $p > \frac{2N}{N+2}$, or,*

*(b) $p = \frac{2N}{N+2}$ and $\alpha > 1/2$,*

*then $\{f^\varepsilon\}$ is precompact in $H^{-1}_{loc}(I\!R^N)$.*

**Proof**. We assume that $\Omega$ is included within the $N$-box, $\mathcal{C}_0 = [-2^{k_0}, 2^{k_0}]^N$. We will consider an orthonormal wavelet basis for $L^2(\Omega)$, $\{\psi_{jk}\}$. This basis may be built using a (finite) wavelet set, $\Psi = \{\psi\}$, supported in $\mathcal{C}_0$, which we will require to belong to $H^1(I\!R^N)$ (consult [6, §10.1], [10, §3.6] or [23] for a brief overview). Specifically, the wavelet basis consists of

$$\psi_{jk}(x) := 2^{kN/2}\psi(2^k x - j), \qquad k \in Z_0^+ := Z^+ - k_0,\ j \in Z^N,\ \psi \in \Psi,$$

which are supported in the dyadic cubes $\mathcal{C}_{jk} := 2^{-k}(\mathcal{C}_0 + j)$; of course, $diam(\mathcal{C}_{jk}) \sim R_k = 2^{-k}$ for all $j$'s, and we consider the wavelet expansion of each $f^\varepsilon$:

$$f^\varepsilon = \sum_{\psi\in\Psi}\sum_{k\in Z_0^+}\sum_{j\in Z^N} \hat{f}^\varepsilon_{jk}\psi_{jk}, \qquad \hat{f}^\varepsilon_{jk} = \int_{\mathcal{C}_{jk}} f^\varepsilon \psi_{jk} dx.$$

The $\psi_{jk}$'s are $H^{-1}$-orthogonal, each of which does not exceed $\|\psi_{jk}\|^2_{H^{-1}} \leq \min\left\{2^{-2k}\int |\widehat{\psi}(\eta)|^2/|\eta|^2, 1\right\}$, and hence

$$\|f^\varepsilon\|^2_{H^{-1}} = \sum_{\psi\in\Psi}\sum_{(j,k)\in(Z^N,Z_0^+)} |\hat{f}^\varepsilon_{jk}|^2 \|\psi_{jk}\|^2_{H^{-1}} \leq Const \sum_{k\in Z_0^+} 2^{-2k} \sum_{j\in Z^N} |\hat{f}^\varepsilon_{jk}|^2.$$

Next we observe that $\cup_j \mathcal{C}_{jk}$ is a covering of disjoint cubes, each of volume of $R_k^N = 2^{-kN}$. Hence, application of (2.4) for $f^\varepsilon \in \vee^{p2,\alpha}$ (with $R = R_k = 2^{-k}$) yields

$$\begin{aligned} \sum_{j\in Z^N} |\hat{f}^\varepsilon_{jk}|^2 &\leq 2^{kN} \sum_{j\in Z^N} \Big( \int_{\mathcal{C}_{jk}} |f^\varepsilon(x)| dx \Big)^2 \qquad (3.16) \\ &\leq Const \cdot 2^{kN} \|f^\varepsilon\|^2_{\vee^{p2,\alpha}} \cdot 2^{-2kN/p'} |1 + k_+|^{-2\alpha}. \end{aligned}$$

It follows that the $f^\varepsilon$'s are bounded in $H^{-1}$. Indeed, Using (3.16) we find the upper-bound

$$\|f^\varepsilon\|^2_{H^{-1}} \leq Const \sum_{k\in Z_0^+} 2^{-2k}\cdot 2^{kN}2^{-2kN/p'}|1+k_+|^{-2\alpha}$$

which shows that $f^\varepsilon$ are $H^{-1}$-bounded if either $(a)$ or $(b)$ holds. Moreover, we have $H^{-1}$-compactness of $\{f^\varepsilon\}$ in view of the *uniform* summability

$$\|\sum_{k>K}\sum_{j\in Z^N}\hat{f}^\varepsilon_{jk}\psi_{jk}\|^2_{H^{-1}} \leq Const. \sum_{k>K} 2^{k(N-2N/p'-2)}|1+k_+|^{-2\alpha} \overset{K\to\infty}{\longrightarrow} 0, \quad \text{uniformly in } \varepsilon.$$

The uniform high-frequency decay (in $H^{-1}$) converts weak compactness in $H^{-1}$ into a strong one. ∎

*Remarks.*

1. The compact imbedding stated in Theorem 3.1 is extended to more general families of measures. Arguing along [21, Theorem 4.3] we find

$$\widetilde{\vee}^{p2,\alpha}(\Omega) \overset{\text{comp}}{\hookrightarrow} H^{-1}(\Omega), \qquad \Big(p-\frac{2N}{N+2}\Big)_+ + \Big(\alpha-\frac{1}{2}\Big)_+ > 0.$$

2. The scale of space $\vee^{pq,\alpha}$ enables to make precise the (compact) embeddings in more general Besov spaces $B^s_\eta(L^r(\Omega))$ spaces (measuring $s$-order of smoothness in $L^r_{\text{loc}}(I\!R^N)$ with secondary index $\eta$). The latter is characterized by a bounded wavelet expansion based on a scaled basis of pre-wavelets $\psi_{jk}(x) = 2^{kN/r}\psi(2^kx-j)$. Assume $\psi$ has certain order of smoothness, say, $s_0$, then [8],[9]

$$\|f\|^\eta_{B^s_\eta(L^r(I\!R^N))} \sim \sum_{k\in Z} 2^{ks\eta}\Big(\sum_{j\in Z^N}|\hat{f}_{jk}|^r\Big)^{\eta/r}, \quad -\infty < s < s_0,\ 1 < r < \infty.$$

Arguing as before we arrive at

**Theorem 3.2** *For a bounded $\Omega \subset I\!R^N$ we have*

$$\vee^{pq,\alpha}(\Omega) \overset{comp}{\hookrightarrow} B^s_\eta(L^q(\Omega)), \qquad \begin{cases} \frac{1}{p} < \frac{1}{q'} - \frac{s}{N}, & \alpha \geq 0 \\ \frac{1}{p} = \frac{1}{q'} - \frac{s}{N}, & \alpha > 1/\eta. \end{cases} \tag{3.17}$$

The case $(\eta,q,s) = (2,2,-1)$ corresponds to Theorem 3.1. The limiting case $(\eta,q,s) = (\infty,p',0)$ yields the (non compact) imbedding

$$\vee^{pp'}(\Omega) \subset B_\infty(L^{p'}(\Omega)). \tag{3.18}$$

Equipped with the scale of spaces $V^{p2,\alpha}$ of Theorem 3.1, we return to examine the gap mentioned earlier concerning $H^{-1}$ compactness. For the full gap of $p$'s, $p \in [\frac{N+2}{2N}, \frac{N}{2}]$, to be $H^{-1}$-compact requires $N/2 < 2$, where we are left with precisely the two relevant cases of two- and three-dimensional problems. We distinguish between the two borderline cases.

• In the two-dimensional case we find that

$$\widetilde{\mathbb{V}}_c^{12,\alpha}(\mathbb{R}^2) \overset{\text{comp}}{\hookrightarrow} H^{-1}_{\text{loc}}(\mathbb{R}^2), \quad \alpha > \frac{1}{2}. \tag{3.19}$$

We recall that for $\alpha > 1$, $\widetilde{M}^{1,\alpha}_{\text{loc}}$ is $H^{-1}_{loc}(\mathbb{R}^2)$-compact, while $\widetilde{M}^{1,1/2}(\mathbb{R}^2) \cap \mathcal{BM}^+_c(\mathbb{R}^2)$ is not. Using (3.19), we are now able to address the open issue of compact imbedding of $\widetilde{M}^{1,\alpha}_{\text{loc}} = \widetilde{\mathbb{V}}^{1\infty,\alpha}$ in the gap $1/2 \le \alpha \le 1$. We conclude that half of this gap, quantified in terms of $\widetilde{\mathbb{V}}^{1q,\alpha}(\mathbb{R}^2)$, $\alpha > \frac{1}{2}$ with the secondary index $1 \le q \le 2$ is $H^{-1}_{\text{loc}}(\mathbb{R}^2)$-compact, and, as we shall see below, the conclusion (3.19) is sharp in the sense that $H^{-1}_{\text{loc}}(\mathbb{R}^2)$-compactness is lost for $\widetilde{\mathbb{V}}^{1q,\frac{1}{2}}(\mathbb{R}^2)$ with the secondary index in the other half, $2 \le q \le \infty$. In particular, we identify as borderline case for $H^{-1}$-compactness, the space $\widetilde{\mathbb{V}}^{12}(\log \widetilde{\mathbb{V}})^{1/2}(\mathbb{R}^2)$ which consists of all measures such that

$$\widetilde{\mathbb{V}}^{12}(\log \widetilde{\mathbb{V}})^{1/2}(\Omega) = \Big\{\mu \mid \sup_{\{B_j\}\subset\mathcal{B}(\Omega)} \sum_j |\log R_j|(|\mu|(B_j))^2 \le Const.\Big\}, \quad \Omega \subset \mathbb{R}^2. \tag{3.20}$$

• In the three-dimensional case we find

$$\widetilde{\mathbb{V}}_c^{p2}(\mathbb{R}^3) \overset{\text{comp}}{\hookrightarrow} H^{-1}_{\text{loc}}(\mathbb{R}^3), \quad p > \frac{6}{5}. \tag{3.21}$$

We recall the different scales of $H^{-1}$-compactness: for Moerry spaces, $\widetilde{M}^p_{\text{loc}}(\mathbb{R}^3) \overset{\text{comp}}{\hookrightarrow} H^{-1}$ for $p > \frac{3}{2}$ while for Lorentz spaces, $L^{p\infty}_{loc}(\mathbb{R}^3) \overset{\text{comp}}{\hookrightarrow} H^{-1}$ for $p > \frac{6}{5}$. Using our new scale of spaces, we can now address the issue of $H^{-1}$-compactness of $\widetilde{M}^p_{loc} = \widetilde{\mathbb{V}}^{p\infty}$ in the gap $3/2 > p > 6/5$. We conclude that for $p > 6/5$, half of this gap, quantified in terms of $\widetilde{\mathbb{V}}^{pq}_{loc}(\mathbb{R}^3)$ with $1 \le q \le 2$, is $H^{-1}$-compact. In particular, we realize that as in the case of Lorentz scale, $p = \frac{6}{5}$ is the 'correct' critical index for $H^{-1}(\mathbb{R}^3)$-compactness, and we identify as borderline case the space $\widetilde{\mathbb{V}}^{\frac{6}{5}2}(\mathbb{R}^3)$ which consists of all $\mu$'s such that

$$\widetilde{\mathbb{V}}^{\frac{6}{5}2}(\Omega) = \Big\{\mu \mid \sup_{\{B_j\}\subset\mathcal{B}(\Omega)} \sum_j \frac{1}{R_j}(|\mu|(B_j))^2 \le Const.\Big\}, \quad \Omega \subset \mathbb{R}^3. \tag{3.22}$$

It is instructive at this point to compare the regularity statement of $\widetilde{\mathbb{V}}^{\frac{6}{5}2}(\mathbb{R}^3)$ vs. the regularity of the 3D borderline case in Morrey scale, $\widetilde{M}^{3/2}(\mathbb{R}^3)$. The latter consists of those $\mu$'s whose total mass over arbitrary balls decays at least linearly with the radius,

$$\widetilde{M}^{3/2}(\Omega) = \Big\{\mu \mid |\mu|(B_R) \leq Const.R, \ \Big\}, \quad \Omega \subset I\!\!R^3.$$

The $\widetilde{\mathbb{V}}^{\frac{6}{5}2}(I\!\!R^3)(\Omega)$-bound in (3.22) allows for a slower decay of the total mass — up to order one-half for a single ball, yet this slower rate should take into account a *collection* of disjoint balls. In general, therefore, the two different bounds are not comparable unless additional information regarding the *asymptotic behavior* of covering balls in (3.22) is provided. For example, an $\widetilde{M}^{3/2}(\Omega)$ bound of $\mu$ yields

$$\|\mu\|^2_{\widetilde{\mathbb{V}}^{\frac{6}{5}2}(\Omega)} \leq \sup_{R_j \leq R_0} \sum_j \frac{1}{R_j}(|\mu|(B_j))^2 \leq \sum_j R_j \cdot \|\mu\|^2_{\widetilde{M}^{3/2}(\Omega)}, \tag{3.23}$$

and hence, if $\Omega$ has a finite *packing measure*, $\pi^{h_1}(\Omega)$, so that it can be packed by covering balls with the finite sum of diameters, we conclude

**Corollary 3.1** *Assume that* $\Omega \subset I\!\!R^3$ *has a finite packing measure,* $\pi^{h_1}(\Omega) < \infty,\ h_1(t) = t$. *Then*

$$\widetilde{M}^{3/2}(\Omega) \subset \widetilde{\mathbb{V}}^{\frac{6}{5}2}(\Omega).$$

# 4 Approximate solutions of Euler's equations

We are concerned with flows of incompressible ideal fluid modeled by the Euler equations

$$\begin{cases} u_t + u \cdot \nabla u = -\nabla p \\ \text{div } u = 0 \\ \text{initial and boundary data,} \end{cases} \tag{4.24}$$

where $u := (u_1, \ldots, u_N)$ and $p$ are the velocity and pressure of the flow. One way to address the question of existence of (weak) solutions for (4.24) is by producing a family of *approximate solutions*, $\{u^\varepsilon(\cdot, t)\}$ and justifying the passage to the limit, say $\varepsilon \downarrow 0$. We recall the definition of *approximate solutions* over any fixed time interval $[0, T]$. We seek a family of incompressible velocity fields, $\{u^\varepsilon\}$, $\text{div}u^\varepsilon = 0$, uniformly bounded in $L^\infty([0,T]; L^2_{\text{loc}}(I\!\!R^N)) \cap Lip((0,T); H^{-L}_{\text{loc}}(I\!\!R^N))$ such that they satisfy the approximate consistency with (4.24). Namely, for any test vector field $\Phi \in C^\infty_c([0,T) \times I\!\!R^N)$ with $\text{div } \Phi = 0$ we have

$$\int_0^T \int_{I\!\!R^N} \Phi_t \cdot u^\varepsilon + (D\Phi\ u^\varepsilon) \cdot u^\varepsilon\ dxdt + \int_{I\!\!R^N} \Phi(x,0) \cdot u^\varepsilon(x,0)\ dx \longrightarrow 0 \text{ as } \varepsilon \to 0. \tag{4.25}$$

The uniform bound in $L^\infty([0,T]; L^2_{\text{loc}}(I\!\!R^N))$ states the uniform bound on the kinetic energy. In the generic case, these weak formulations hold in some negative Sobolev space tested against vector fields in $H^s_c([0,T) \times I\!\!R^N)$. Together with the $L^2$-energybound, it follows that $u^\varepsilon$ has the Lip regularity with a uniform bound in $Lip((0,T);H^{-L}_{\rm loc}(I\!R^N))$ for some $L=L(s,N)>1$, e.g., [11],[18]. This (weak) regularity in time enables us to define the manner in which $u^\varepsilon$ assumes prescribed initial data.

The $L^2$-energybound implies that we can extract a weak-* converging subsequence, $\{u^{\varepsilon_k}\}\rightharpoonup u$ in $L^\infty([0,T];L^2_{\rm loc}(I\!R^N))$, and thus we are facing one of two possibilities. Either there is *strong* $L^2$-convergence $u^{\varepsilon_k}(\cdot,t)\to u(\cdot,t)$ in $L^1[0,T]$, so that by passing to the limit (in both the linear and quadratic terms) in (4.25), $u(\cdot,t)$ is found to be a weak solution of (4.24). The other possibility is lack of strong convergence, $u^\varepsilon\rightharpoonup u$. In this case the $L^2$ energy concentrates on a subset $E\subset\Omega\times[0,T]$ characterized by a positive *reduced defect measure* introduced in $\theta(E)>0$, [12],

$$\theta(E):=\limsup\int_E |u(x,t)-u^\varepsilon(x,t)|^2dxdt. \tag{4.26}$$

Outside this concentration set $\text{limsup}_{\varepsilon\to 0}\int_{E^c}|u^\varepsilon-u|^2dxdt=0$. Greengard & Thomann [17] have shown that the concentration set $E$ has Hausdorff dimension $H(E)\geq 1$. Upper bounds on the 2D Hausdoff dimension $H(E)$ can be found in [26].

The phenomena of energy concentration does not exclude the possibility of convergence to a weak solution. DiPerna & Majda initiated in [11],[12],[13] the study of the *concentration-cancelation* phenomena, where subtle cancellation justify the passage to limit $u_i^{\varepsilon_k}u_j^{\varepsilon_k}\rightharpoonup u_iu_j$, $i\neq j$, so that despite the concentration of energy, weak-$*\lim u^{\varepsilon_k}=u$ is found to be a weak solution of (4.24).

It is physically relevant to classify many approximate flows into one of the two scenarios outlined above according to the behavior of their vorticity fields, $\boldsymbol{\omega}^\varepsilon(\cdot,t):=\nabla\times u^\varepsilon(\cdot,t)$. A sharp criterion for strong $L^2$-convergence was introduced in the recent work [21]. The so called $H^{-1}$-*stability* criterion requires the associated vorticity field $\boldsymbol{\omega}^\varepsilon(\cdot,t)$ to form a precompact subset in $C((0,T),H^{-1}_{\rm loc}(I\!R^N))$. The main result [21, Theorem 1.1] states that an $H^{-1}$-stable family of approximate solutions, $\{u^\varepsilon\}$, admits a subsequence which is strongly convergent to a weak solution in $L^\infty([0,T],L^2_{\rm loc}(I\!R^N))$.

We will utilize the $H^{-1}$ stability criterion to study the strong convergence of approximate Euler solutions. In particular, our new refined scale of spaces, $\vee^{pq,\alpha}(I\!R^N)$, will enable us to 'approach' the borderline cases which separate the phenomena of concentration-cancelation. We distinguish between two-dimensional and three-dimensional flows.

## 4.1 The 2D Euler equations

Incompressible flows in two space dimensions become considerably simpler (than the $N>2$-case), since the 2D vorticity equation is reduced to the scalar transport equation

$$\omega_t+u\cdot\nabla\omega=0. \tag{4.27}$$

It is governed by a divergence-free velocity field, $u$, which is recovered by the Biot-Savart law $u=K*\omega$ with $K(\xi):=\xi^\perp/(2\pi|\xi|^2)$. It follows that any *rearrangement invariant* space, $X$, is a regularity space for the vortcity equation (4.27), so that $X$-regularity of $\omega^\varepsilon(\cdot,t)$ is retained in time. Thus, consider a specific example family of approximate Euler solutions, $\{u^\varepsilon(\cdot,t)\}$, associated with the mollified initial data, $u_0^\varepsilon = K_\varepsilon * \omega_0$, where $K_\varepsilon$ denotes the mollified kernel $K_\varepsilon := \eta_\varepsilon * K$. It follows – consult [21, Corollary 2.2] for the precise details, that if the initial vorticity, $\omega_0$, belongs to such rearrangement-invariant space $X$, $X_{\text{loc}} \overset{\text{comp}}{\hookrightarrow} H^{-1}_{\text{loc}}(I\!\!R^2)$, then $H^{-1}$-stability is retained for later times, and hence $\{u^\varepsilon\}$ has a strong limit, $u(\cdot,t)$, which is a weak solution associated with the initial velocity $u_0 = K * \omega_0 \in X$ without concentrations.

The 2D rearrangement-invariant examples of $H^{-1}$-compactness revisited in [21] (generalizing [3],[24],[20]), include

{i} Orlicz spaces, $L(\log L)^\alpha_c(I\!\!R^2), \alpha > 1/2$; and the slightly larger

{ii} Lorentz spaces $L_c^{(1,q)}(I\!\!R^2), q < 2$.

We also mention the borderline cases which are not compactly imbedded in $H^{-1}_{loc}(I\!\!R^2)$,

{iii} $X = L(\log L)^{1/2}_c(I\!\!R^2)$ and $X = L_c^{(1,2)}(I\!\!R^2)$

Despite the lack of compactness in these borderline cases, it was shown in [21, Theorem 2.2 & Theorem 2.4] that special $X$-sequences of approximate vorticities corresponding to mollified initial data in these borderline cases, $\omega_0^\varepsilon = \eta_\varepsilon * \omega_0$, $\omega_0 \in X$, are in fact $H^{-1}_{loc}(I\!\!R^2)$-compact.

The 2D problem beyond rearrangement-invariant spaces was studied in [21, §3] in terms of Morrey spaces, $M_c^{1,\alpha}(I\!\!R^2)$, which are compactly imbedded in $H^{-1}_{loc}(I\!\!R^2)$ for $\alpha > 1$. The study of Morrey spaces in this context was motivated by the DiPerna-Majda conjecture on the concentration-cacellation phenomenon of *one-signed* vorticities. Majada, [22], has shown how the Morrey regularity in $\widetilde{M}_c^{1,1/2}(I\!\!R^2)$ of such one-signed vorticities plays a fundamental role in his simplified proof of the concentration-cancellation argument of Delort [14]. The new ladder of spaces, $\vee^{1q,\alpha}(I\!\!R^2)$, provides us with a more precise information on the regularity of one-signed measures which could not be classified in terms of the missing gap in the ladder of Morrey spaces, $M^{1,\alpha}(I\!\!R^2)$, $1/2 < \alpha < 1$.

We begin with an immediate consequence of our main Theorem 3.1 regarding approximate vorticities, $\omega^\varepsilon(\cdot,t) \in X_\alpha := \widetilde{V}_c^{(12,\alpha)}(I\!\!R^2)$. Taking into account the definition of approximate solutions, we have that $\{\omega^\varepsilon\}$ are uniformly bounded,

$$\{\omega^\varepsilon\} \hookrightarrow Lip((0,T), H^{-L-1}_{\text{loc}}(I\!\!R^2)) \cap L^\infty((0,T), X_\alpha)),$$

where according to (3.19), $X_\alpha \overset{\text{comp}}{\hookrightarrow} H^{-1}_{\text{loc}}(I\!\!R^2) \overset{\text{comp}}{\hookrightarrow} H^{-L-1}$. It follows that $\{\omega^\varepsilon\} \overset{\text{comp}}{\hookrightarrow} C((0,T), H^{-1}_{\text{loc}}(I\!\!R^2))$ and by our $H^{-1}$-stability result [21], we conclude

**Corollary 4.1** *Let $\{u^\varepsilon\}$ be a family of approximate solutions of the 2D Euler equations (4.24), and assume that the corresponding sequence of vorticities $\{\omega^\varepsilon\}$ is uniformly bounded in $L^\infty([0,T]; \widetilde{\vee}_c^{(12,\alpha)}(I\!\!R^2))$, with $\alpha > 1/2$. Then $\{u^\varepsilon\}$ is strongly compact in $L^\infty([0,T]; L^2_{loc}(I\!\!R^2))$, and has a strong limit, $u(\cdot,t)$, which is a weak solution with no concentrations.*

Seeking a strategy for obtaining apriori $\vee^{12,\alpha}$-bounds of the type required in the last corollary, we are led to the following.

*Question.* Consider a sequence of approximate vorticities, $\omega^\varepsilon(\cdot,t)$, corresponding to mollified initial data, $\omega_0^\varepsilon = \eta_\varepsilon * \omega_0$ with $\omega_0 \in \widetilde{\vee}_c^{12,\alpha}(I\!R^2)$. Does the sequence $\{\omega^\varepsilon(\cdot,t)\}$ remain in $\widetilde{\vee}_c^{12,\alpha}(I\!R^2)$ for $t>0$?

Though the general question remains open, we offer one possible strategy for obtaining apriori $\vee$-type bounds in the special case of one-signed vorticities. To this end, we let $H(\omega)$ denote the "pseudo-energy"

$$H(\omega) := -\frac{1}{2\pi}\int\int_{I\!R^2\times I\!R^2} \log|x-y|\omega(x)\omega(y)dxdy,$$

noting that it is an invariant quantity associated with smooth vorticities,

$$H(\omega(t)) = H(\omega_0).$$

Indeed, expressed in terms of the streamfunction, $\psi = 1/2\pi \log|x| * \omega$, the velocity is given by $u = \nabla^\perp\psi$, and the energy associated with the 2D flow reads

$$\int_{B_R(0)} |u|^2 dx = \int_{B_R(0)} \nabla^\perp\psi\cdot\nabla^\perp\psi dx = -\int_{B_R(0)} \omega\psi dx + \int_{\partial B_R(0)} \nabla^\perp\psi\cdot \mathbf{t}\psi ds,$$

and hence, assuming a far-field behavior which is invariant in time (– there is no far-field decay of this boundary term), we conclude that in fact $H(\omega(\cdot,t))$ measures the invariance of the total energy $\int |u(\cdot,t)|^2dx$.

Equipped with the invariance of pseudo-energy we now turn to consider the $\widetilde{\vee}^{12,\alpha}$-bound of one-signed vorticities.

**Lemma 4.1** *Let $\{u^\varepsilon\}$ be a family of approximate solutions of the 2D Euler equations with one signed measured vorticities, $\{\omega_0^\varepsilon \in \mathcal{BM}_c^+\}$. Then*

$$\|\omega^\varepsilon(\cdot,t)\|_{\widetilde{\vee}^{12}(\log\widetilde{\vee})_c^{1/2}(I\!R^2)} \le Const. \tag{4.28}$$

**Proof**. We consider an arbitrary collection of disjoint balls, $\{B_j\}_j$ with sufficiently small radii, $B_j = B_{R_j}(x_j)$, $R_j < 1/2$. We partition the energy between its self-induced part, $H_{si}$, and the interaction energy, $H_{ie}$, [4]

$$\begin{aligned} H(\omega^\varepsilon(\cdot,t)) &= -\frac{1}{2\pi}\sum_k \int\int_{B_{R_j}\times B_{R_j}} \log|x-y|d\omega^\varepsilon(x,t)d\omega^\varepsilon(y,t) + \\ & \quad -\frac{1}{2\pi}\sum_{j\ne k}\int\int_{B_{R_j}\times B_{R_k}} \log|x-y|d\omega^\varepsilon(x,t)d\omega^\varepsilon(y,t) =: \\ &=: H_{si}(\omega(t)) + H_{ie}(\omega(t)). \end{aligned}$$

First we note a lower bound on the interaction energy either when $\omega^\varepsilon(\cdot,t)$ remains compactly supported, say in $B_{R_t}(0)$, so that $\log|x-y| \le (\log|2R_t|)_+$, or, following [22], using the fact that $(\log|x-y|)_+ \le 2(|x|^2+|y|^2)$. In either case we find $H_{ie}(\omega^\varepsilon(\cdot,t))$ to be bounded from below; for example, in the second case we find

$$
\begin{aligned}
-H_{ie}(\omega^\varepsilon(\cdot,t)) &\le \frac{1}{\pi}\int\int_{I\!R^2\times I\!R^2}(|x|^2+|y|^2)d\omega^\varepsilon(x,t)d\omega^\varepsilon(y,t)\\
&\le \frac{2}{\pi}I_0(\omega^\varepsilon(\cdot,t)\cdot I_2(\omega^\varepsilon(\cdot,t) \le Const_0.
\end{aligned}
$$

The last uniform bound follows from the fact that the first two moments, $I_0(\omega^\varepsilon(\cdot,t)$ and $I_2(\omega^\varepsilon(\cdot,t)$, are global invariants of 2D flows (or at least bounded quantities for approximate flows).

Second, we note that the $\widetilde{\vee}^{12}(\log\widetilde{\vee})^{1/2}$-bound of $\omega^\varepsilon$ is a *lower bound* for the self-induced energy: indeed, in view of the positivity of $\omega^\varepsilon$,

$$
\begin{aligned}
&\frac{1}{2\pi}\sum_j |\log 2R_j|\Big(\int_{B_j}|d\omega^\varepsilon(x,t)|\Big)^2 \le\\
&\qquad -\frac{1}{2\pi}\sum_j\int\int_{B_{R_j}\times B_{R_j}}\log|x-y|d\omega^\varepsilon(x)d\omega^\varepsilon(y) = H_{si}(\omega^\varepsilon(\cdot,t)).
\end{aligned}
$$

The $\vee$-bound (4.28) follows from the last two estimates,

$$
\begin{aligned}
\frac{1}{2\pi}\|\omega^\varepsilon(\cdot,t)\|^2_{\widetilde{\vee}^{12}(\log\widetilde{\vee})^{1/2}(I\!R^2)} &\le H_{si}(\omega^\varepsilon(\cdot,t)) = H(\omega^\varepsilon(\cdot,t)) - H_{ie}(\omega^\varepsilon(\cdot,t)) \le\\
&\le H(\omega_0) + Const_0.
\end{aligned}
$$

∎

According to Theorem 3.1, $\widetilde{\vee}^{12,\alpha}(I\!R^2)$ are compactly imbedded in $H^{-1}_{loc}(I\!R^2)$ for $\alpha > 1/2$, and by the main stability result of [21], therefore, no concentration phenomenon occurs in this range when $\|\omega^\varepsilon(\cdot,t)\|_{\widetilde{\vee}^{12,\alpha}(I\!R^2)} \le Const.$ $\alpha > 1/2$. In particular, the $\widetilde{\vee}^{12}(\log\widetilde{\vee})^{1/2}$ regularity of one-signed measures is in fact a borderline case and, analogous with our previous discussion of the borderline cases of Orlicz and Lorentz cases, we raise the following.

*Question.* Consider a sequence of approximate vorticities, corresponding to mollified initial data, $\omega_0^\varepsilon = \eta_\varepsilon * \omega_0$, $\omega_0 \in \widetilde{\vee}^{12}(\log\widetilde{\vee})^{1/2}_c(I\!R^2)$. Does the sequence $\{\omega^\varepsilon(\cdot,t)\}$ remain compact in $H^{-1}_{loc}(I\!R^2)$ for $t>0$?

An affirmative answer to this question implies that for 2D initial vorticities with one-signed (in fact – more general) measures, one can construct a solution by a limiting argument which avoids the phenomena of concentration (consult [21, Lemma 2.3] regarding the issue of temporal continuity).

As we noted before in the context of the borderline cases $X = L(\log L)^{1/2}_c(\mathbb{R}^2)$ and $X = L^{(1,2)}_c(\mathbb{R}^2)$, they both lack $H^{-1}_{\rm loc}(\mathbb{R}^2)$-compactness, hence only special $X$-sequences are expected to form compact subsets in $H^{-1}_{\rm loc}(\mathbb{R}^2)$, e.g., approximate vorticities corresponding to mollified initial data. Similarly, we note that only special $\vee^{12}(\log\vee)^{1/2}(\mathbb{R}^2)$-sequence can be expected to form $H^{-1}$-compact sequences. The following counterexample due to DiPerna & Majda [13, Proposition 3.1], demonstrates a family of steady vorticities, $\{\omega^\varepsilon\}$, which are positive and hence uniformly bounded in $\widetilde{\vee}^{12}(\log\widetilde{\vee})^{1/2}$, yet it lacks $H^{-1}$-compactness. To this end, pick a non-negative $C_0^\infty(0,1)$ radial vorticity, $\omega(r)$, and consider its dilations

$$\omega^\varepsilon(x) := \frac{1}{\varepsilon^2\sqrt{|\log\varepsilon|}}\omega\Big(\frac{|x|}{\varepsilon}\Big), \qquad \Gamma(r) := \int_0^r s\omega(s) < \infty.$$

A straightforward computations shows the induced velocity field satisfies the steady Euler equations

$$u^\varepsilon(x) = \frac{1}{\varepsilon\sqrt{|\log\varepsilon|}}u(\frac{x}{\varepsilon}), \quad u(x) = \frac{x^\perp}{|x|^2}\Gamma(|x|),$$

with finite kinetic energy, and for which $u_i^\varepsilon(x)u_j^\varepsilon(x) \rightharpoonup \pi\Gamma^2(\infty)\delta(x)\delta_{ij}$.

The lack of $H^{-1}_{\rm loc}(\mathbb{R}^2)$-compactness for *general* sequences in the borderline case $X = \widetilde{\vee}^{12}(\log\widetilde{\vee})^{1/2}_c(\mathbb{R}^2)$ indicates the possibility of energy concentration, and in this context we show that if energy concentration does take place than the $\widetilde{\vee}^{12}(\log\widetilde{\vee})^{1/2}$ bound is sufficient to guarantee the concentration-cancellation phenomena. The following is a generalization of Delort's result [14].

**Theorem 4.1** *Let $\{u^\varepsilon(\cdot,t)\}$ be a family of approximate solutions of the 2D Euler equations (4.24), and assume that the corresponding sequence of vorticities $\{\omega^\varepsilon\}$ is uniformly bounded in $L^\infty([0,T];\widetilde{\vee}^{12}(\log\widetilde{\vee})^\alpha_c(\mathbb{R}^2))$, $\alpha > 0$. Then the $L^2$-weak limit, $u^\varepsilon \rightharpoonup u(\cdot,t)$ is a finite-energy solution of the 2D Euler equation (4.24).*

**Proof**. A weak formulation of the 2D Euler's equations (4.27)

$$\omega_t + K * \omega \cdot \nabla\omega = 0,$$

reads, consult [27],

$$\int_0^\infty\int_{\mathbb{R}^2}\phi_t\omega^\varepsilon(x,t)dxdt + \int_0^\infty\int_{\mathbb{R}^2\times\mathbb{R}^2}H_\phi(x,y,t)\omega^\varepsilon(x,t)\omega^\varepsilon(y,t)dxdydt + \\ + \int_{\mathbb{R}^2}\phi(x,0)\omega_0^\varepsilon(x)dx = 0, \qquad \forall\phi\in C_c^\infty([0,\infty)\times\mathbb{R}^2),$$

where the kernel $H_\phi(x,y,t)$ is given by

$$H_\phi(x,y,t) := \frac{\nabla\phi(x,t)-\nabla\phi(y,t)}{4\pi|x-y|}\cdot\frac{(x-y)^\perp}{|x-y|}.$$

By a density argument we may restrict our attention to test function of the form $\phi(x,t) = \psi(t)\varphi(x)$. We let $\rho(|x|) \in C_0^\infty(0,2)$ be a positive cut-off function with $\rho(|x|) \equiv 1$ for $|x| \leq 1$. The main issue is passage to a limit in the quadratic term (corresponding to the mixed term $weak - * \lim u_1^\varepsilon u_2^\varepsilon$), which is decomposed in the by-now standard fashion, consult [14],[22, §2], [27],...

$$
\begin{aligned}
&\int_0^\infty \int_{I\!R^2\times I\!R^2} H_\phi(x,y,t)\omega^\varepsilon(x,t)\omega^\varepsilon(y,t)dxdydt = \\
&\quad = \int_0^\infty \int_{I\!R^2\times I\!R^2} \psi(t)\Big(1-\rho\Big(\frac{|x-y|}{\delta}\Big)\Big) H_\varphi(x,y)\omega^\varepsilon(x,t)\omega^\varepsilon(y,t)dxdydt + \\
&\quad + \int_0^\infty \int_{I\!R^2\times I\!R^2} \psi(t)\rho\Big(\frac{|x-y|}{\delta}\Big) H_\varphi(x,y)\omega^\varepsilon(x,t)\omega^\varepsilon(y,t)dxdydt =: \\
&\quad =: I_\delta(\omega^\varepsilon) + J_\delta(\omega^\varepsilon).
\end{aligned}
$$

By Delort's lemma, [14, Proposition 1.2.3], $\psi(t)\Big(1-\rho\Big(\frac{|x-y|}{\delta}\Big)\Big)H_\varphi(x,y)$ is a 'nice' kernel such that

$$
\lim_{\varepsilon\downarrow 0} I_\delta(\omega^\varepsilon) = \int_0^\infty \int_{I\!R^2\times I\!R^2} \psi(t)\Big(1-\rho\Big(\frac{|x-y|}{\delta}\Big)\Big) H_\varphi(x,y)d\omega(x,t)d\omega(y,t)dt.
$$

It remains to estimate the behavior of $J_\delta(\omega^\varepsilon)$ which is supported near the singularity along the diagonal $x = y$, and it is here that the $\widetilde{\vee}^{12,\alpha}$-bound plays an essential role. To this end, we cover $I\!R^2$ with a net of $2\delta \times 2\delta$ cubes, $\mathcal{C}_j = 2\delta\mathcal{C}(\cdot + 2\delta j),\ j \in Z^2$, with $\mathcal{C}$ denoting the 2D unit cube. Decomposing the integration in $J_\delta(\omega^\varepsilon)$ over $I\!R^2 = \cup_j \mathcal{C}_j$, we find

$$
\begin{aligned}
J_\delta(\omega^\varepsilon) &= \int_0^\infty \psi(t) \sum_{j,k} \int_{(x,y)\in(\mathcal{C}_j\times\mathcal{C}_k)} \rho\Big(\frac{|x-y|}{\delta}\Big) H_\varphi(x,y)\omega^\varepsilon(x,t)\omega^\varepsilon(y,t)dxdydt \leq \\
&\leq C_\varphi \cdot \int_0^\infty |\psi(t)| \sum_{j,k} \int_{\substack{(x,y)\in(\mathcal{C}_j\times\mathcal{C}_k)\\ |x-y|\leq 2\delta}} |\omega^\varepsilon(x,t)|\cdot|\omega^\varepsilon(y,t)|dxdydt, \quad C_\varphi := \|H_\varphi\|_{L^\infty}.
\end{aligned}
$$

For each cell $\mathcal{C}_j$, only its immediate neighboring cells, $\mathcal{C}_k,\ |k-j|_\infty \leq 1$, participate in the summation on the right of (4.29), so that $|x-y| \leq 2\delta$ whenever $(x,y) \in (\mathcal{C}_j, \mathcal{C}_k)$. For each $j$ (respectively $k$) there are precisely nine such neighboring cells (including the cell $k = j$ itself) which contribute to the self-induced energy. Since $\rho H_\varphi$ is bounded along the diagonal we find, in view of (2.4)

$$
\begin{aligned}
J_\delta(\omega^\varepsilon) &\leq C_\varphi \cdot \int_0^\infty |\psi(t)| \Big[\sum_{|j-k|\leq 1} \frac{1}{2}\Big(\int_{x\in\mathcal{C}_j} |\omega^\varepsilon(x,t)|dx\Big)^2 + \sum_{|j-k|\leq 1} \frac{1}{2}\Big(\int_{x\in\mathcal{C}_k} |\omega^\varepsilon(y,t)|dy\Big)^2\Big)\Big] dt \leq \\
&\leq 9\cdot C_\varphi \cdot \int_0^\infty |\psi(t)| \Big(\int_{x\in\mathcal{C}_j} |\omega^\varepsilon(x,t)|dx\Big)^2 dt \leq \\
&\leq 9C_\varphi \|\psi\|_{L^\infty} \cdot \|\omega^\varepsilon(x,t)\|^2_{L^1_{\text{loc}}(I\!R_t;\widetilde{\vee}^{12}(\log\widetilde{\vee})^\alpha)} \times |\log\rho|^{-2\alpha}, \qquad C_\varphi = \|H_\varphi\|_{L^\infty}.
\end{aligned}
$$

It follows that $J_\delta(\omega^\varepsilon)$ tends to zero uniformly in $\varepsilon$, $|J_\delta(\omega^\varepsilon)| \leq Const|\log\rho|^{-2\alpha} \underset{\rho\to 0}{\longrightarrow} 0$, and we conclude that $H_\phi(x,y,t)\omega^\varepsilon(x,t)\omega^\varepsilon(y,t) \rightharpoonup H_\phi(x,y,t)d\omega(x,t)d\omega(y,t)$. ∎

## 4.2 The 3D Euler equations

According to the compact imbedding (3.15), a family of 3D vorticities, $\{\boldsymbol{\omega}^\varepsilon(\cdot,t)\}$, which is uniformly bounded in $L^\infty([0,T];\widetilde{M}^p_{loc}(\mathbb{R}^3))$ with $p > 3/2$, induces a velocity field with $L^2$-strong limit, $u(\cdot,t)$, which is a global weak solution of the 3D Euler equations, consult [21, Theorem 4.5]. Remark that unlike the 2D problem, however, Morrey space estimates do not have the physical interpretation as circulation decay estimates. And moreover, there is no known strategy of obtaining apriori estimates on the $\widetilde{M}^p$-size of the vorticity, $\|\boldsymbol{\omega}(\cdot,t)\|_{\widetilde{M}^p(\mathbb{R}^3)}$, in time. We want to show that our new scale of spaces offers a better tool to handle the issue of compactness in terms of physically relevant invariant quantities.

As in the 2D case, we begin with the following.

**Corollary 4.2** *Let $\{u^\varepsilon\}$ be a family of approximate solutions of the three-dimensional Euler equations (4.24), and assume that the corresponding sequence of vorticities $\{\boldsymbol{\omega}^\varepsilon\}$ is uniformly bounded in $L^\infty([0,T];(\widetilde{\vee}^{p2}(\mathbb{R}^3))$, with $p > 6/5$. Then, $\{u^\varepsilon\}$ is strongly compact in $L^\infty([0,T];L^2_{loc}(\mathbb{R}^3))$, and hence it has a strong limit, $u(\cdot,t)$, which is a weak solution with no concentrations.*

There is no known strategy to obtain apriori $\vee^{p2}(\mathbb{R}^3)$-bounds on the vorticity, and there is no a priori reason to expect that they are invariants of 3D flows. There is one notable exception, however, which is linked precisely to the borderline case of $\vee^{p2}(\mathbb{R}^3)$ with $p = 6/5$. We explore this exceptional case below. First we recall the one special 3D invariant which is the pseudo-energy (the Coulomb energy)

$$H(\boldsymbol{\omega}(x,t)) := \frac{1}{8\pi}\int\int_{\mathbb{R}^3\times\mathbb{R}^3} \frac{\langle\boldsymbol{\omega}(x,t),\boldsymbol{\omega}(y,t)\rangle}{|x-y|}dxdy = H(\boldsymbol{\omega}_0).$$

Next, we cover space with a 3D lattice, $\mathbb{R}^3 = \cup_j \mathcal{C}_j$, and as before, we partition the energy, $H(\boldsymbol{\omega}(x,t)) = H_{si}(\boldsymbol{\omega}(x,t)) + H_{ie}(\boldsymbol{\omega}(x,t))$, into its self induced, short range part, $H_{si}(\boldsymbol{\omega}(x,t))$, and long range interaction energy, $H_{ie}(\boldsymbol{\omega}(x,t))$, namely

$$\begin{aligned}
H_{si}(\boldsymbol{\omega}(x,t)) &= \frac{1}{8\pi}\sum_j \int\int_{\mathcal{C}_j\times\mathcal{C}_j} \frac{\langle\boldsymbol{\omega}(x,t),\boldsymbol{\omega}(y,t)\rangle}{|x-y|}dxdy \\
H_{ie}(\boldsymbol{\omega}(x,t)) &= \frac{1}{8\pi}\sum_{j\neq k} \int\int_{\mathcal{C}_j\times\mathcal{C}_k} \frac{\langle\boldsymbol{\omega}(x,t),\boldsymbol{\omega}(y,t)\rangle}{|x-y|}dxdy.
\end{aligned}$$

To proceed, we make two claims regarding *lower bounds* of the two portions of the pseudo-energy, similar to the 2D configuration (but much harder to prove):

(i) A lower bound on the interaction energy

$$H_{ie}(\boldsymbol{\omega}^\varepsilon(x,t)) = \frac{1}{8\pi}\sum_{j\neq k}\int\int_{\mathcal{C}_j\times\mathcal{C}_k}\frac{\langle \boldsymbol{\omega}^\varepsilon(x,t), \boldsymbol{\omega}^\varepsilon(y,t)\rangle}{|x-y|}dxdy \geq -Const_{ie}. \tag{4.29}$$

(ii) For sufficiently small cubes, $\mathcal{C}_j$,

$$\boldsymbol{\omega}^\varepsilon(x,t) \sim ⨍_{\mathcal{C}_j}\boldsymbol{\omega}^\varepsilon(y,t)dy, \qquad x\in\mathcal{C}_j \tag{4.30}$$

which leads to a lower bound of the self-induced energy

$$\begin{aligned} H_{si}(\omega^\varepsilon(\cdot,t)) &= \frac{1}{8\pi}\sum_j\int\int_{\mathcal{C}_j\times\mathcal{C}_j}\frac{\langle \boldsymbol{\omega}(x,t), \boldsymbol{\omega}(y,t)\rangle}{|x-y|}dxdy \geq \\ &\geq \frac{1}{8\pi}\sum_j\frac{1}{2R_j}\Big(\int_{\mathcal{C}_j}|\boldsymbol{\omega}^\varepsilon(x,t)|dx\Big)^2. \end{aligned}$$

The last two estimates yield the $V^{\frac{6}{5}2}(I\!\!R^3)$-bound on the vorticity, consult (3.22),

$$\begin{aligned} &\frac{1}{16\pi}\sum_j\frac{1}{R_j}\Big(\int_{\mathcal{C}_j}|\boldsymbol{\omega}^\varepsilon(x,t)|dx\Big)^2 \leq \\ &\qquad \leq H_{si}(\boldsymbol{\omega}^\varepsilon(\cdot,t)) = H(\boldsymbol{\omega}^\varepsilon(\cdot,t)) - H_{ie}(\boldsymbol{\omega}^\varepsilon(\cdot,t)) \leq H_0 + Const_{ie}, \quad H_0 = H(\boldsymbol{\omega}_0^\varepsilon). \end{aligned}$$

The new ladder of spaces is establishing a direct linkage between question related to the global configuration of the 3D pseudo-energy and the borderline case of $\vee^{\frac{6}{5}2}(I\!\!R^3)$-regularity. Similar to the 2D framework we are now led to inquire about the $H^{-1}(I\!\!R^3)$ compactness of this borderline case.

*Question.* Consider a sequence of approximate vorticities, $\boldsymbol{\omega}^\varepsilon(\cdot,t)\in L^\infty([0,T],\vee^{\frac{6}{5}2}(I\!\!R^3))$. What are the possible configurations of the pseudo-energy so that the sequence $\{\boldsymbol{\omega}^\varepsilon(\cdot,t)\}$, is compact in $L^\infty([0,T],H^{-1}_{\rm loc}(I\!\!R^3))$?

We note that an answer to this question maps a possible strategy of constructing solutions to 3D Euler equation in the large. The estimates claimed in (4.29), (4.30) demonstrated this issue. For a detailed discussion on the configurations of the self-induced energy the interaction energy and the relation to vortex stretching we refer to Chorin [4, Chapter 5].

We conclude this section by pointing out one such strategy which leads to the desired $\vee^{\frac{6}{5}2}$-bound in the 3D case. To this end we let

$$\boldsymbol{\xi}(x,t) := \frac{\boldsymbol{\omega}^\varepsilon(x,t)}{|\boldsymbol{\omega}^\varepsilon(x,t)|}, \quad \boldsymbol{\omega}^\varepsilon(x,t)\neq 0,$$

denote the direction of the vorticity $\boldsymbol{\omega}^\varepsilon$. The stretching effect of $\boldsymbol{\omega}^\varepsilon$ is controlled by the difference $|\boldsymbol{\xi}^\varepsilon(x,t)-\boldsymbol{\xi}^\varepsilon(y,t)|$ and we make

**Assumption 4.1** *There exist constants $\delta>0$ and $\theta=\theta_\delta<1$, such that whenever $|\boldsymbol{\omega}^\varepsilon(x,t)|,|\boldsymbol{\omega}^\varepsilon(y,t)|>K_0$, there holds*

$$|\boldsymbol{\xi}^\varepsilon(x,t)-\boldsymbol{\xi}^\varepsilon(y,t)|\leq\sqrt{2}\theta,\qquad \forall|x-y|\leq\delta. \tag{4.31}$$

Squaring (4.31) yields $2-2\langle\boldsymbol{\xi}^\varepsilon(x,t),\boldsymbol{\xi}^\varepsilon(y,t)\rangle=|\boldsymbol{\xi}^\varepsilon(x,t)-\boldsymbol{\xi}^\varepsilon(y,t)|^2\leq 2\theta^2$, and hence, whenever $|\boldsymbol{\omega}^\varepsilon(x,t)|,|\boldsymbol{\omega}^\varepsilon(y,t)|>K_0$, we have

$$\langle\boldsymbol{\omega}^\varepsilon(x,t),\boldsymbol{\omega}^\varepsilon(y,t)\rangle\geq(1-\theta^2)\cdot|\boldsymbol{\omega}^\varepsilon(x,t)|\cdot|\boldsymbol{\omega}^\varepsilon(y,t)|,\qquad |x-y|\leq\delta. \tag{4.32}$$

Thus, under assumption (4.31) there is a local alignment of the *direction of the vorticity*, $\boldsymbol{\xi}^\varepsilon(\cdot,t)$, whenever its magnitude, $|\boldsymbol{\omega}^\varepsilon(\cdot,t)|$, becomes too large. Assumption 4.1 is inspired by Constantin & Fefferman, [5], who proved the existence of 3D Navier-Stokes solutions under the short-range alignment assumption

$$|\boldsymbol{\xi}^\varepsilon(x,t)-\boldsymbol{\xi}^\varepsilon(y,t)|\leq\frac{|x-y|}{\delta},\qquad |x-y|\leq\delta,\ \ |\boldsymbol{\omega}^\varepsilon(x,t)|,|\boldsymbol{\omega}^\varepsilon(y,t)|>K_0.$$

Equipped with the alignment assumption 4.1 we prove that $\boldsymbol{\omega}^\varepsilon(\cdot,t)$ remains uniformly bounded in the borderline space $X_3=\vee_c^{\frac{6}{5}2}(I\!R^3)$

**Theorem 4.2** *Let $\{u^\varepsilon(\cdot,t)\}$ be a family of approximate solutions of the 3D Euler equations (4.24). Assume that the corresponding sequence of vorticities, $\{\boldsymbol{\omega}^\varepsilon(\cdot,t)\}$, is compactly supported and satisfies the local alignment condition (4.31). Then the following holds,*

$$\|\boldsymbol{\omega}^\varepsilon(\cdot,t)\|_{\vee^{\frac{6}{5}2}(\Omega)}\leq Const_T,\qquad \Omega\subset I\!R^3,\ \ t\leq T. \tag{4.33}$$

*Remark.* The requirement of $\boldsymbol{\omega}^\varepsilon(\cdot,t)$ having compact support is made for simplicity and could be replaced by a weaker requirement of fast enough decay at infinity.

**Proof**. We begin by partitioning the total energy between its short-range, self-induced part, and its long-range interaction energy, $H(\boldsymbol{\omega}^\varepsilon(\cdot,t))=H_{si}(\boldsymbol{\omega}^\varepsilon(\cdot,t))+H_{ie}(\boldsymbol{\omega}^\varepsilon(\cdot,t))$. The partition is taken at a scale level $\delta$ dictated by the alignment assumption in (4.31),

$$H_{si}(\boldsymbol{\omega}^\varepsilon(\cdot,t)):=\frac{1}{8\pi}\int\int_{I\!R^3\times I\!R^3}1_{|x-y|\leq\delta}\frac{\langle\boldsymbol{\omega}^\varepsilon(x,t),\boldsymbol{\omega}^\varepsilon(y,t)\rangle}{|x-y|}dxdy,$$
$$H_{ie}(\boldsymbol{\omega}^\varepsilon(\cdot,t)):=\frac{1}{8\pi}\int\int_{I\!R^3\times I\!R^3}1_{|x-y|\geq\delta}\frac{\langle\boldsymbol{\omega}^\varepsilon(x,t),\boldsymbol{\omega}^\varepsilon(y,t)\rangle}{|x-y|}dxdy,.$$

In the 3D case we have the advantage the interaction energy is lower-bounded (by $H(\boldsymbol{\omega}^\varepsilon(\cdot,t))=-H_0$), or equivalentaly, that $H_{si}(\boldsymbol{\omega}^\varepsilon)\leq 2H_0$. Indeed, computing the 3D

Fourier transform, $\eta(\xi) := \mathcal{F}(|x|^{-1}1_{|x|<\delta}) = |\xi|^{-2}(1-\cos(|\xi|\delta))$, yields for the weighted $L^2_\eta$ norm,

$$\begin{aligned} H_{si}(\boldsymbol{\omega}^\varepsilon(\cdot,t)) &= \frac{1}{8\pi}\|\hat{\boldsymbol{\omega}}^\varepsilon(\xi,t)\|^2_{L^2_{\eta(\xi)}} \leq \\ &\leq \frac{2}{8\pi}\int_{\xi\in I\!R^3}\frac{|\hat{\boldsymbol{\omega}}^\varepsilon(\xi,t)|^2}{|\xi|^2}d\xi = 2H(\boldsymbol{\omega}^\varepsilon(\cdot,t)) \leq 2H_0. \end{aligned} \tag{4.34}$$

Next, we split the vorticity between its bounded and unbounded parts at 'height' $K_0$

$$\boldsymbol{\omega}^\varepsilon(x,t) = \boldsymbol{\omega}^\varepsilon(x,t)1_{\Omega\cap\{x\ |\ |\boldsymbol{\omega}^\varepsilon(x,t)|\leq K\}} + \boldsymbol{\omega}^\varepsilon(x,t)1_{\Omega\cap\{x\ |\ \boldsymbol{\omega}^\varepsilon(x,t)|<K_0\}} =: \boldsymbol{\omega}^\varepsilon_- + \boldsymbol{\omega}^\varepsilon_+,$$

and we show that the bounded part of the vorticity, $\boldsymbol{\omega}^\varepsilon_-(\cdot,t)$, has a finite contribution to the self-induced energy. We start by expanding

$$H_{si}(\boldsymbol{\omega}^\varepsilon(\cdot,t)) \equiv H_{si}(\boldsymbol{\omega}^\varepsilon_+(\cdot,t)) + H_{si}(\boldsymbol{\omega}^\varepsilon_-(\cdot,t)) + 2H_{si}(\boldsymbol{\omega}^\varepsilon_-(\cdot,t),\boldsymbol{\omega}^\varepsilon_+(\cdot,t)),$$

with the third term on the right denoting the bilinear positive form ( — positivity follows along the lines of (4.34) or consult [19, Theorem 9.8])

$$H_{si}(\boldsymbol{f},\boldsymbol{g}) := \frac{1}{8\pi}\int\int_{|x-y|\leq\delta}\frac{\langle \boldsymbol{f}(x),\boldsymbol{g}(y)\rangle}{|x-y|}dxdy.$$

Cauchy-Schwartz inequality then yields

$$2|H_{si}(\boldsymbol{\omega}^\varepsilon_-(\cdot,t),\boldsymbol{\omega}^\varepsilon_+(\cdot,t))| \leq \frac{1}{2}H_{si}(\boldsymbol{\omega}^\varepsilon_+(\cdot,t)) + 8H_{si}(\boldsymbol{\omega}^\varepsilon_-(\cdot,t)),$$

and in view of (4.34) we end up with the upper-bound

$$\begin{aligned} \frac{1}{2}H_{si}(\boldsymbol{\omega}^\varepsilon_+(\cdot,t)) &\leq H_{si}(\boldsymbol{\omega}^\varepsilon(\cdot,t)) + 7H_{si}(\boldsymbol{\omega}^\varepsilon_-(\cdot,t)) \leq \\ &\leq 2H_0 + 7Const_{K_0}, \qquad H_0 = H(\boldsymbol{\omega}^\varepsilon(\cdot,0)). \end{aligned} \tag{4.35}$$

Here we used the fact that $\boldsymbol{\omega}^\varepsilon(\cdot,t)$ are compactly supported so that

$$H_{si}(\boldsymbol{\omega}^\varepsilon_-(\cdot,t)) \leq \frac{K_0^2}{8\pi}\int_{x\in supp\ \boldsymbol{\omega}^\varepsilon(\cdot,t)}\int_{y\in B_\rho(x)}\frac{1}{|x-y|}dydx \leq Const_{K_0},$$

with, say, $Const_{K_0} \sim K_0^2|diam(supp\ \boldsymbol{\omega}^\varepsilon(\cdot,t))|^3\delta^2$. Of course, one can relax the requirement of compact support, asking a fast enough decay of $\boldsymbol{\omega}^\varepsilon(x,t),\ |x|\to\infty$.

Given a collection of disjoint balls, $\{B_j = B_{R_j}(x_j)\}$, we claim the short range part of the energy controls the $\vee$-size of $\boldsymbol{\omega}^\varepsilon_+$, when measured over all balls with radii $R_j \leq R_0 < \delta/4$; indeed, in view of (4.32), our alignment assumption implies

$$\begin{aligned} H_{si}(\boldsymbol{\omega}_+^\varepsilon(\cdot,t)) &\geq \frac{1}{8\pi}\sum\int\int_{(x,y)\in B_j\times B_j}\frac{(1-\theta^2)|\boldsymbol{\omega}_+^\varepsilon(x,t)|\cdot|\boldsymbol{\omega}_+^\varepsilon(y,t)|}{|x-y|}dxdy \geq \\ &\geq \frac{1-\theta^2}{16\pi}\sum_j\frac{1}{R_j}\Big(\int_{B_j}|\boldsymbol{\omega}_+^\varepsilon(x,t)|dx\Big)^2, \end{aligned}$$

and by varying over all collections of such balls we find a lower bound for the self-induced part of the energy in terms of its $\vee^{\frac{6}{5}2}$-norm

$$H_{si}(\boldsymbol{\omega}_+^\varepsilon(\cdot,t)) \geq \frac{1-\theta^2}{16\pi}\|\boldsymbol{\omega}_+^\varepsilon\|^2_{\vee^{\frac{6}{5}2}(\Omega)}.$$

Using this estimate together with (4.35), the asserted $\vee^{\frac{6}{5}2}$-bound follows,

$$\|\boldsymbol{\omega}_+^\varepsilon\|^2_{\vee^{\frac{6}{5}2}(\Omega)} \leq \frac{32\pi}{1-\theta^2}(2H_0 + 7Const_{K_0}).$$

∎

The $\vee^{\frac{6}{5}2}_{\rm loc}(I\!R^3)$-bound derived in Theorem 4.2 implies that $\{\boldsymbol{\omega}^\varepsilon(\cdot,t)\}$ is uniformly bounded in $H^{-1}_{\rm loc}(I\!R^3)$. This in turn can be strengthened into $H^{-1}_{\rm loc}$-compactness, for example, as long as the velocity field remains uniformly $L^{p>2}_{\rm loc}$-bounded. The proof is essentially an application of Murat Lemma, [25]. Arguing along the lines of [21, Theorem 4.6] we conclude

**Corollary 4.3** *Let $\{u^\varepsilon(\cdot,t)\}$ be a family of approximate solutions of the 3D Euler equations (4.24) such that $\{u^\varepsilon\}$ is uniformly bounded in $L^\infty((0,T),L^p(I\!R^3))$ with $p>2$, and assume that the compactly supported vorticities, $\{\boldsymbol{\omega}^\varepsilon(\cdot,t)\}$ satisfy the local alignment condition (4.31). Then $\{u^\varepsilon\}$ is strongly compact in $L^\infty((0,T),L^2_{loc}(I\!R^3))$, and hence it has a strong limit, $u(\cdot,t)$, which is a weak solution of (4.24).*

We close by noting that there is no known strategy to guarantee the $L^p_{\rm loc}(I\!R^3)$-bound on the velocity for $p>2$.

*Added to the proofs.* If we let $Pf(x) := \sum_j ⨍_{B_j}|f|\cdot\chi_{B_j}(x)$ denote the Haar projection of $f$ subject to the partition $\{B_j\}$, then a straightforward computation shows $\|f\|_{V^{pp}} = sup\{\|Pf\|_{L^p(\Omega)} \mid \{B_j\}\subset\mathcal{B}(\Omega)\}$, and by a density argument therefore, $V^{pp}\subset L^p$. It follows that $V^{pq}$ forms the scale of interpolation spaces between $V^{pp}=L^p$ and $M^p$.